\documentclass[12 pt]{amsart}
\usepackage{amsmath}
\usepackage{amssymb}
\usepackage{amscd}
\usepackage{enumerate}
\usepackage{mathtools}
\usepackage{dashrule}
\usepackage{bm}
\usepackage{mathrsfs}
\usepackage{tabu}
\usepackage{tikz}
\usepackage{eucal}
\usepackage[T1]{fontenc}
\usepackage{comment}
\PassOptionsToPackage{hyphens}{url}\usepackage{hyperref}
\usepackage[all,cmtip]{xy}
\usepackage[letterpaper, top = 2cm, bottom = 2cm,right = 4cm, left = 4cm]{geometry}

\theoremstyle{theorem}
\newtheorem{theorem}{Theorem}
\newtheorem{lemma}{Lemma}

\theoremstyle{definition}
\newtheorem{definition}{Definition}
\newtheorem{remark}{Remark}
\newtheorem{example}{Example}

\numberwithin{equation}{section}

\begin{document}
\sloppy


\title{Ziegler Spectrum and Krull Gabriel Dimension}


 \author{Daniel L\'{o}pez-Aguayo}
 \address{Centro de Ciencias Matem\'{a}ticas
UNAM,Morelia,Mexico}
\email{dlopez@matmor.unam.mx}

\maketitle

\begin{abstract}
These notes are based on a talk given at the Summer School ``\emph{Infinite-dimensional representations of finite-dimensional algebras}" held at the University of Manchester in September 2015. They intend to provide a brief introduction to the notion of Ziegler Spectrum and Krull-Gabriel dimension. 
\end{abstract}

\bigskip

Let $k$ be a field and let $A$ be an associative $k$-algebra. Let $\operatorname{Mod}$-A (respectively A-$\operatorname{Mod}$) denote the category of right $A$-modules (respectively left $A$-modules) and let $D$: A-$\operatorname{Mod}$ $\rightarrow$ $\operatorname{Mod}$-A be the functor $D(-)=\operatorname{Hom}_{k}(-,k)$. Throughout these notes $\mathbf{Ab}$ will denote the category of abelian groups.

\begin{definition} A monomorphism $M \stackrel{f} \hookrightarrow N$ in $\operatorname{Mod}$-A is pure if after tensoring with every left $A$-module remains injective. More precisely, if $f \otimes 1_{X}: M \otimes_{A} X \rightarrow N \otimes_{A} X$ is injective for each left $A$-module $X$.
\end{definition}

\begin{definition} A module $M$ is pure injective if every pure monomorphism $f: M \hookrightarrow L$ splits, i.e. there exists an $A$-module map $g: L \rightarrow M$ such that $gf=id_{M}$.
\end{definition}

\begin{definition} Let $\mathcal{E}: 0 \rightarrow  B \stackrel{\alpha} \rightarrow C \stackrel{\beta} \rightarrow D \rightarrow 0$ be a short exact sequence in $\operatorname{Mod}$-A. We say $\mathcal{E}$ is pure-exact if for every left $A$-module $X$ the homomorphism $B \otimes_{A} X \stackrel{\alpha \otimes 1_{X}} \longrightarrow C \otimes_{A} X$ is injective. 
\end{definition}

\begin{remark} If the contravariant functor $\operatorname{Hom}_{A}(-,Q)$ is exact on  pure exact sequences starting at $Q$ then $Q$ is pure injective.
\end{remark}

\begin{example} Duals are always pure-injective: let $M$ be a left $A$-module and let $0 \rightarrow D(M) \rightarrow N \rightarrow L \rightarrow 0$ be pure-exact. Then $0 \rightarrow D(M) \otimes_{A} M \rightarrow N \otimes_{A} M \rightarrow L \otimes_{A} M \rightarrow 0$ is exact. Applying the contravariant exact functor $\operatorname{Hom}_{k}(-,k)$ yields the following short exact sequence: \\

\begin{center}
$0 \rightarrow \operatorname{Hom}_{k}(L \otimes_{A} M,k) \rightarrow \operatorname{Hom}_{k}(N \otimes_{A} M,k) \rightarrow \operatorname{Hom}_{k}(D(M) \otimes_{A} M,k) \rightarrow 0$
\end{center}

\smallskip

Recall that $\operatorname{Hom}_{A}(-,D(M))  \cong \operatorname{Hom}_{k}( - \otimes_{A} M,k)$. Thus we get a short exact sequence: 

\smallskip

\begin{center}
$0 \rightarrow \operatorname{Hom}_{A}(L,D(M)) \rightarrow \operatorname{Hom}_{A}(N,D(M)) \rightarrow \operatorname{Hom}_{A}(D(M),D(M)) \rightarrow 0$
\end{center}

It follows that $D(M)$ is pure exact. In particular, any finite dimensional module is pure injective (as $M \cong D^{2}(M))$. 
\end{example}

\begin{definition} An indecomposable $A$-module of infinite dimension is said to be \textbf{generic} provided $M$ considered as an $\operatorname{End}_{A}(M)$-module is of finite length (i.e endofinite).
\end{definition}

\begin{example} An example of a generic module is the following representation of the Kronecker quiver:
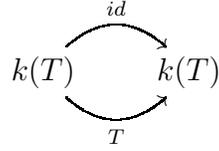
\begin{displaymath}
\xymatrix{
k(T) \ar@/^1.5pc/[r]^{id} \ar@/_1.5pc/[r]_{T}  &  k(T)
}
\end{displaymath}

here, $k(T)$ denotes the field of rational functions in one variable $T$ over the base field $k$.

\end{example}

\begin{example} Any endofinite module is pure-injective. In particular, generic modules are pure-injective.
\end{example}

\begin{example} The evaluation map $M \rightarrow M^{\ast \ast}$ is a pure monomorphism. See \cite[Corollary 1.3.16]{Prest1}.
\end{example}

\begin{definition}
The Ziegler spectrum of an algebra $A$, $Zg A$, is defined to be the set of isoclasses of indecomposable pure-injective $A$-modules. Later, we will see that this set can be endowed with a topology.
\end{definition}

\begin{example} Up to isomorphism, the indecomposable pure injective $\mathbb{Z}$-modules \cite{Kaplansky} are: 

\begin{enumerate}[(a)]
\item finite groups $\mathbb{Z}/{p^n\mathbb{Z}}$,
\item Pr\"ufer groups $\mathbb{Z}_{p^{\infty}}$,
\item the $p$-adic integers $\overline{\mathbb{Z}_{(p)}}$,
\item $\mathbb{Q}$
\end{enumerate}

where $p$ ranges over primes and $n$ over positive integers. This determines completely $Zg \mathbb{Z}$. Remark that the fact that $(c)$ and $(d)$ are pure injective follows from the well known facts that these are divisible abelian groups (recall that a $\mathbb{Z}$-module is injective iff it is divisible) thus they are injective and therefore pure injective.
\end{example}

\begin{remark} Some facts about the Ziegler spectrum of a finite dimensional algebra (see \cite{Krause}):
\begin{enumerate}[(a)]
\item $\{M\}$ is open if and only if $M$ is finite dimensional. This can be shown using Auslander-Reiten theory.
\item The finite dimensional points form a dense subset.
\item If $R$ is of finite representation type then the $ZgR$ is discrete (i.e. has the discrete topology). [If $R$ is of finite representation type then by \cite[Proposition 4.5.22]{Prest1} it follows that every indecomposable module is endofinite  and thus by \cite[Proposition 5.1.12]{Prest1} such a module is a closed point of $ZgR$. Since there are finitely many points in $ZgR$, the claim follows.]
\end{enumerate}
\end{remark}

Let $\mathscr{A}$ be an additive category. Given an object $X \in \mathscr{A}$ we denote by $H_{X}$ the representable functor:

\begin{center}
$H_{X}:=Hom_{\mathscr{A}}(X,-): \mathscr{A} \rightarrow \mathbf{Ab}$, $Y \mapsto Hom_{\mathscr{A}}(X,Y)$.
\end{center}

\begin{definition} A functor $F: \mathscr{A} \rightarrow \mathbf{Ab}$ is finitely presented if there exists an exact sequence:

\begin{center}
$H_{Y} \rightarrow H_{X} \rightarrow F \rightarrow 0$
\end{center}

in the functor category $(\mathscr{A},\mathbf{Ab})$.
\end{definition}

\begin{definition} An $A$-module $M$ is finitely presented if there is an exact sequence:

\begin{center}
$A^{m} \rightarrow A^{n} \rightarrow M \rightarrow 0$
\end{center}

where $m,n \in \mathbb{N}$. 
\end{definition}

\begin{theorem}\cite[6.1]{Auslander} Let $M$ be a right $R$-module. The functor $M \otimes_{R} -: R-\operatorname{mod} \rightarrow \mathbf{Ab}$ is finitely presented if and only if $M$ is finitely presented.
\end{theorem}

\begin{definition} An additive functor $F$: Mod-A $\rightarrow \mathbf{Ab}$ is said to be coherent if it commutes with direct limits and products. The coherent functors form an abelian category whose morphisms are the natural transformations.
\end{definition}

In the next list, modules are taken in the category $\operatorname{Mod}$-A  where A is a finite dimensional $k$- algebra. The following list provides examples of coherent functors \cite{William}.
\begin{enumerate} 
\item The representable functor $\operatorname{Hom}_{A}(M,-)$ where $M$ is finite dimensional.
\item The tensor product $M \otimes_{A} -$ where $M$ is a finite-dimensional right $A$-module.
\item The derived functor $Ext^{i}(X,-)$ where $X$ is finite dimensional.
\item The functor $Tor_{i}(N,-)$ where $N$ is finite dimensional.
\end{enumerate} 

\smallskip

We now state a theorem which gives a way of putting a topology on $Zg A$ via coherent functors. We won't assume compact spaces are Hausdorff; in fact the Ziegler spectrum is seldom Hausdorff \cite[8.2.12]{Prest1}

\begin{theorem} (Ziegler) The sets $\{ M \in Zg A: F(M) \neq 0\}$ with $F$ a coherent functor, form a base of open sets for a topology on $Zg A$. With this topology $Zg A$ becomes a compact topological space.
\end{theorem}

\begin{center}
\textbf{Serre subcategories and localization}
\end{center}

\begin{definition} Let $\mathcal{A}$ be an abelian category. Let $\mathcal{S}$ be a full subcategory of $\mathcal{A}$. We say that $S$ is a Serre subcategory (or a Serre class) of $\mathcal{A}$ such that whenever $0 \rightarrow A \rightarrow B \rightarrow C \rightarrow 0$ is a short exact sequence in $\mathcal{A}$ then $B \in \mathcal{S}$ if and only if $A,C \in \mathcal{S}$.
\end{definition}

\begin{example}
All torsion groups, all finitely generated groups, all finite groups, all $p$-groups in the category $\mathbf{Ab}$ are Serre subcategories.
\end{example}

\begin{example}
The subclass of $\operatorname{Mod}-A$ consisting of Noetherian (resp. Artinian) right $A$-modules is a Serre class.
\end{example}

\begin{example} If $F: \mathcal{A} \rightarrow \mathcal{B}$ is a functor, then the kernel of $F$ is defined to be the full subcategory with objects:

\begin{center}
$\operatorname{Ker}(F)=\{A \in ob(\mathcal{A}): F(A)=0\}$
\end{center}

If $F$ is exact, then the kernel of $F$ is a Serre class. This is immediate from the fact that if a sequence $0 \rightarrow F(X) \rightarrow 0$ is exact if and only if $F(X)=0$ i.e iff $X \in ker(F)$.
\end{example}

\begin{example} Check that the class $\{M \in \operatorname{R-Mod}: I^{n}M=0 \ \textrm{for some} \ n\geq 1\}$, for a fixed ideal $I \subseteq R$, is Serre.
\end{example}

Let $\mathscr{A}$ be an abelian category.

\begin{definition} A subobject of an object $X$ is an object $Y$ together with a monomorphism $i: Y \hookrightarrow X$.
\end{definition}

\begin{definition} A quotient object of $Y$ is an object $Z$ with an epimorphism $p: Y \rightarrow Z$.
\end{definition}

\begin{definition} A subquotient object of $Y$ is a quotient object of a subobject of $Y$.
\end{definition}

\begin{definition} For a subobject $X \subset Y$ define the quotient object $Z:=Y/X$ to be the cokernel of a monomorphism $f: X \rightarrow Y$. We write $Y \subseteq X$.
\end{definition}

The notion of a Serre class $\mathcal{S}$ of $\mathscr{A}$ will be used to define a \textbf{quotient category} in the sense of Serre-Grothendieck. We sketch the construction. See \cite[pp.498-505]{Faith} for more details.

Let $A,B$ be objects of $\mathscr{A}$ and let $A',B'$ be subobjects of $A$ and $B$ respectively (i.e. $A' \subseteq A$, $B' \subseteq B$). Let $I$ denote the class of ordered pairs $\langle A',B' \rangle$ where $A' \subseteq A$, $B' \subseteq B$ and $A/A',B' \in \mathcal{S}$. Order $I$ by: \\

\begin{center}
$\langle A',B' \rangle \leq \langle A'',B'' \rangle$ iff $A' \supseteq A''$ and $B' \subseteq B''$.
\end{center}

\smallskip

One can verify that $I$ is a directed set. 

\smallskip

Note that if $\langle A',B' \rangle \leq \langle A'',B'' \rangle$ then there is an induced map: \\

\begin{center}
$\operatorname{Hom}_{\mathscr{A}}(A',B/B') \rightarrow \operatorname{Hom}_{\mathscr{A}}(A'',B/B') \rightarrow \operatorname{Hom}_{\mathscr{A}}(A'',B/B'')$
\end{center}

\smallskip

Then $\{\operatorname{Hom}_{\mathscr{A}}(A',B/B'): A/A' \in \mathcal{S},B' \in \mathcal{S}\}$ is a directed system of abelian groups indexed by $I$ so it has a direct limit. \\

Let $\mathscr{A}$ be an abelian category and let $\mathcal{S}$ be a Serre subcategory. The \textbf{quotient category} $\mathscr{A}/\mathcal{S}$ is defined as follows. The objects of $\mathscr{A}/\mathcal{S}$ are the same objects of $\mathscr{A}$ and the morphisms are defined as:

\begin{center}
$\operatorname{Hom}_{\mathscr{A}/\mathcal{S}}(A,B):=\varinjlim (A',B/B')$
\end{center}

where the limit is over subobjects $A' \subseteq A$, $B' \subseteq B$ such that $A/A',B' \in \mathcal{S}$.

There is clearly a \textbf{quotient functor} $\pi: \mathscr{A} \rightarrow \mathscr{A}/\mathcal{S}$: we define $\pi(A)=A$ for every object $A \in \mathscr{A}$ and let $\pi(f: M \rightarrow N)$ be the image of $f$ in the direct system that defines $\operatorname{Hom}_{\mathscr{A}/\mathcal{S}}(M,N)$.

\begin{example}
The quotient of $\mathbf{Ab}$ modulo its Serre subcategory of torsion groups is the category of $\mathbf{Q}$-vector spaces. 
\end{example}

An alternative way of constructing the quotient category with respect a Serre class $\mathcal{S}$ is to use localization where the multiplicative system is given by the set $C_{\mathcal{S}}=\{f \in \operatorname{Mor}(\mathscr{A}): \operatorname{ker}(f), \operatorname{coker}(f) \in Ob(\mathcal{S})\}$. Then we let: 

\begin{center}
$\mathscr{A}/\mathcal{S}:=\mathscr{A}[C_{\mathcal{S}}^{-1}]$
\end{center}

\smallskip

\begin{definition} Let $\mathcal{S}$ be a Serre subcategory. We say $\mathcal{S}$ is a \textbf{localizing subcategory} of $\mathscr{A}$, if the quotient functor $T: \mathscr{A} \rightarrow \mathscr{A}/\mathcal{S}$ has a right adjoint $G: \mathscr{A}/\mathcal{S} \rightarrow \mathscr{A}$; the functor $G$ is called a section functor. This means that for every two objects $A \in \mathscr{A}, B \in \mathscr{A}/\mathcal{S}$ there are bijections:

\begin{center}
$\xymatrix{
\operatorname{Hom}_{\mathscr{A}/\mathcal{S}}(T(A),B) \ar@<-4pt>[r]_{\phi_{A,B}} & \operatorname{Hom}_{\mathscr{A}}(A,G(B)) \ar@<-4pt>[l]_{\psi_{A,B}} 
}$
\end{center}
natural in $A$ and $B$. The functor $T$ is also called the \textbf{localization functor}.
\end{definition}

The following theorem ensures that $\mathscr{A}/\mathcal{S}$ is an abelian category, so the quotient category is not the same as factoring by an ideal of $\mathcal{S}$ (the process used to obtain the stable module category $\underline{\operatorname{mod}} \ \Lambda$).

\begin{theorem} Let $\mathscr{A}$ be an abelian category and let $\mathcal{S}$ be a Serre subcategory of $\mathscr{A}$. The quotient category $\mathscr{A}/\mathcal{S}$ is abelian and the localization functor $F: \mathscr{A} \rightarrow \mathscr{A}/\mathcal{S}$ is exact. The kernel of $F=\{A \in \mathscr{A}: F(A)=0\}$ is $\mathcal{S}$. Every exact functor $G: \mathscr{A} \rightarrow \mathscr{A'}$ where $\mathscr{A'}$ is abelian category and $\mathcal{S} \subseteq ker(G)$ factors uniquely through $F: \mathscr{A} \rightarrow \mathscr{A}/\mathcal{S}$ via an exact and faithful functor.
\end{theorem}

There is full and faithful embedding of $\operatorname{Mod}-R$ into the functor category $(R-\operatorname{mod},\mathbf{Ab})$ where $R-\operatorname{mod}$ is the category of finitely presented left $R$-modules which is given on objects by taking the $R$-module $M_{R}$ to the tensor product $M \otimes_{R} - : R-mod \rightarrow \mathbf{Ab}$ (see \cite[Theorem 2.2]{Prest3}). 

\begin{definition} Let $\mathscr{C}$ be a preadditive category. A generating set for $\mathscr{C}$ is a collection $\mathcal{A}$ of objects of $\mathscr{C}$ such that for every nonzero morphism $f: B \rightarrow C$ in $\mathscr{C}$ there is $G \in \mathcal{A}$ and a morphism $g: G \rightarrow B$ such that $fg \neq 0$. If $\mathcal{A}=\{G\}$ then we say that $G$ is a generator.
\end{definition}

\begin{example} In the category $\mathbf{Ab}$ the abelian group $\mathbb{Z}$ is a generator. In the category $\mathbf{Set}$ any singleton is a generator.
\end{example}

\begin{definition} A Grothendieck category is an abelian category which has arbitrary coproducts, in which direct limits are exact and which has a generator.
\end{definition}

The following are examples of Grothendieck categories.
\begin{enumerate}
\item Given any associative unital ring $R$, the category $\operatorname{Mod}-R$ is a Grothendieck category. In particular, $\mathbf{Ab}$ is a Grothendieck category. 
\item Given a small category $\mathscr{C}$ and a Grothendieck category $\mathscr{A}$, the functor category $Funct(\mathscr{C},\mathscr{A})$ is a Grothendieck category. If $\mathscr{C}$ is preadditive, then the functor category $\operatorname{Add}(\mathscr{C},\mathscr{A})$ of all additive functors is a Grothendieck category as well.
\item If $\mathscr{A}$ is a Grothendieck category and $\mathcal{S}$ is a localizing subcategory of $\mathscr{A}$, then the Serre quotient category $\mathscr{A}/\mathcal{S}$ is a Grothendieck category.
\end{enumerate}

\begin{definition} Let $\mathscr{C}$ be an additive category with direct limits. An object $A$ of a category $\mathscr{C}$ is \textbf{finitely presented} if the representable functor $(A,-):=\operatorname{Hom}_{\mathscr{C}}(A,-)$ commutes with direct limits, meaning that if $((C_{\lambda})_{\lambda},(f_{\lambda \mu}: C_{\lambda} \rightarrow C_{\mu})_{\lambda<\mu})$ is any directed system in $\mathscr{C}$ and if $(C,(f_{\lambda}: C_{\lambda} \rightarrow C)_{\lambda})$ is the direct limit of this system then the direct limit of the directed system $((A,C_{\lambda}))_{\lambda}, (A,f_{\lambda \mu}: (A,C_{\lambda}) \rightarrow (A,C_{\mu}))_{\lambda<\mu})$ in $\mathbf{Set}$ is $((A,C),(A,f_{\lambda})_{\lambda})$.
\end{definition}

 Let $\mathscr{A}$ be a Grothendieck category. A full subcategory $\mathcal{T}$ is a $\textbf{torsion subcategory}$ if it is closed under extensions, quotient objects and arbitrary coproducts. If it is also closed under subobjects then it is a \textbf{hereditary} torsion subcategory (or subclass). The objects of $\mathcal{T}$ are referred to as \textbf{torsion} and those in the corresponding torsionfree subcategory $\mathcal{F}=\{D: (\mathcal{T},D)=0\}=\{D: \operatorname{Hom}_{\mathscr{A}}(T,D)=0 \ \text{for all T in $\mathcal{T}$} \}$ are the \textbf{torsionfree} objects. It is then the case that $\mathcal{T}=\{C: (C,\mathcal{F})=0\}$ and the pair $\tau=(\mathcal{T},\mathcal{F})$ is referred to as a \textbf{torsion theory}, which is said to be \textbf{hereditary} if $\mathcal{T}$ is closed under subobjects.

\begin{definition} Let $\tau=(\mathcal{T},\mathcal{F})$ be a (hereditary) torsion theory on $\mathscr{A}$. We say $\tau$ is of \textbf{finite type} if the torsionfree class $\mathcal{F}$ is closed under direct limits. Equivalently, if the torsion class $\mathcal{T}$ is generated as such by the class of finitely presented torsion objects.
\end{definition}

If $\mathscr{A}$ is a functor category, then the closure under direct limits, $\overrightarrow{\mathcal{S}}$, of any Serre subcategory $\mathcal{S}$ of $\mathscr{A}^{fp}$ is a hereditary torsion class in $\mathscr{A}$, that is, a class closed under subobjects, factor objects, extensions and arbitrary coproducts. 

Hereditary torsion theories, $\tau$, of finite type on the functor category $(mod-\mathscr{A},\mathbf{Ab})$ correspond bijectively to the Serre subcategories, $\mathcal{S}$, of the finitely presented functor category $\operatorname{fun}-\mathscr{A}=(mod-\mathscr{A},\mathbf{Ab})^{fp}$, the maps being $\tau \mapsto \mathcal{T}_{\tau} \cap fun-\mathscr{A}$ and the inverse map sends $S \mapsto \overrightarrow{S}$.

\begin{definition} A Grothendieck category $\mathscr{A}$ is said to  be \textbf{locally coherent} provided that $\mathscr{A}$ has a generating set of finitely presented objects and the full subcategory of finitely presented objects in $\mathscr{A}$ is abelian.
\end{definition}

\begin{definition} Let $X$ be an object in an abelian category $\mathscr{C}$. We say that $X$ has finite length if there exists a filtration:

\begin{center}
$0=X_{0} \subset X_{1} \subset ... \subset X_{n-1} \subset X_{n}=X$
\end{center}

such that $X_{i}/X_{i-1}$ is simple for all $i$. Such a filtration is called a Jordan-H\"older series of $X$.
\end{definition}

Let $\mathscr{A}$ be a locally coherent Grothendieck category. Let $\mathcal{S}$ be the Serre subcategory of all finitely presented objects of $\mathscr{A}$ of finite length. Then $\overrightarrow{\mathcal{S}}$ is a finite type torsion class in $\mathscr{A}$ and we may form the localization $\mathscr{A}_{1}=\mathscr{A}/\overrightarrow{\mathcal{S}}$. By \cite[Corollary 3.6]{Prest3} $\mathscr{A}_{1}$ is again a locally coherent category, so we can repeat the process with $\mathscr{A}_{1}$ in place of $\mathscr{A}$, and obtain $\mathscr{A}_{}$, etc, transfinitely. We obtain a sequence of locally coherent categories $\mathscr{A}_{\alpha}$ indexed by ordinals. 
If $\alpha$ is the smallest ordinal such that $\mathscr{A}_{\alpha}=0$ then we say that the \textbf{Krull-Gabriel dimension} of $\mathscr{A}$ is $\alpha$ and write $KGdim\mathscr{A}=\alpha$. If eventually we reach a category with no finitely presented simple objects then we set $KGdim(\mathscr{A})=\infty$. \\

If $\mathscr{A}$ is the functor category $(mod-R, \mathbf{Ab})$ then we write $KG(R)$ for $KG dim(\operatorname{mod}-R,\mathbf{Ab})$. There is a connection between this dimension and representation type of Artin algebras \cite[p.23]{Prest3}: 

\begin{enumerate}
\item $KG(R)=0$ if and only if $R$ is of finite representation type. (e.g. the truncated polynomial algebra $R=k[x]/(x^{n})$).
\item There is no Artin algebra whose Krull-Gabriel dimension is $1$. 
\item Algebras with Krull-Gabriel dimension $2$ include the tame hereditary algebras (see \cite{Geigle}).
\item Wild-finite dimensional algebras have infinite Krull-Gabriel dimension (see \cite[pp.281-2]{Prest2}); for example, the path algebra of the extended three Kronecker quiver: 

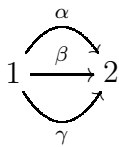
\begin{displaymath}
\xymatrix{
1 \ar@/^1.5pc/[r]^{\alpha} \ar@/_1.5pc/[r]_{\gamma} \ar[r]^{\beta} &  2
}
\end{displaymath}
has infinite Krull-Gabriel dimension.

\item For hereditary Artin algebras (an algebra of global dimension $\leq 1$ that is finitely generated as a module over its Artinian center; e.g. the path algebra of an acyclic finite quiver) the possible representation types -finite,tame/domestic, wild- correspond to the values $0,2,\infty$ for Krull-Gabriel dimension.
\item The domestic string algebras $\Lambda_{n}$ have $KG(\Lambda_{n})=n+1$ (see \cite{Burke}).
\end{enumerate}

\end{document}